\documentclass{article}
\usepackage[utf8]{inputenc}
\usepackage[top=1in, bottom=1in, left=1.25in, right=1.25in]{geometry}
\usepackage{amsmath}
\usepackage{graphicx} 

\usepackage{latexsym}
\usepackage{enumerate}
\usepackage{epsf,epsfig,amsfonts,a4wide}
\usepackage{mathtools}
\usepackage{amsfonts}
\usepackage{url}
\usepackage{amsmath,amssymb,amsthm}
\usepackage{epsf,epsfig,amsfonts,graphicx}
\usepackage{a4wide}

\usepackage[centerlast]{caption}
\usepackage{graphics,amsmath,amssymb,amscd}
\usepackage{amsbsy}
\usepackage{amsthm}
\usepackage{gensymb}
\usepackage{float}
\usepackage{graphicx}
\usepackage[colorlinks=true,linkcolor=black,citecolor=black]{hyperref}
\usepackage{multicol}
\usepackage{color}

\usepackage{epstopdf}

\newcommand{\be}{\begin{equation}}
\newcommand{\ee}{\end{equation}}

\usepackage{amsthm}
\usepackage{mathdots}

\usepackage{amscd}

\title{\textbf{On stability of planar solutions of double averaged restricted elliptic three-body problem} }
\author{Anatoly  Neishtadt$^{1,2}$ \footnote{Corresponding author.  E-mail addresses:  {\it A.Neishtadt@lboro.ac.uk} (A.Neishtadt), {\it K.Sheng@lboro.ac.uk} (K.Sheng), {\it vvsidorenko@list.ru}  (V.Sidorenko)}\, , Kaicheng Sheng$^1$, Vladislav Sidorenko$^{3,4}$
 \\ 
$^1$ Loughborough University, Loughborough, LE11 3TU, UK\\
$^2$ Space Research Institute, Moscow, 117997, Russia\\
$^3$ Keldysh Institute for Applied Mathematics, Moscow, 125047, Russia\\
$^4$ Moscow Institute of Physics and Technology, Moscow Region, 141701, Russia}

\date{\today}

\begin{document}

\maketitle

\begin{abstract}

Double averaged  planar restricted elliptic three-body problem has a two-parametric family of stable equilibria. We show that these equilibria are stable in the linear approximation  as equilibria of the  double averaged spatial restricted elliptic three-body problem. They are Lyapunov stable for all values of parameters but, possibly,  parameters from some finite set of analytic curves.

\end{abstract}

\section{Introduction}

Averaged models play an important role in celestial mechanics. We consider the restricted three-body (a star, a planet, and an asteroid) problem \cite{Sze} when the mass of the planet is much smaller  than the mass of the star. In this case one can use  averaging over motions of the system star - planet and of the asteroid (double averaging).

 In  case of  circular orbits of the star and the planet the double averaged problem is integrable \cite{moiseev}.  This problem is considered in \cite{lidov} under assumption  that the distance between  the asteroid and  the star is much smaller  than the distance between  the planet and the star (Hill's approximation).   The complete analytical study of this problem is given. Results of this study  are rediscovered in \cite{kozai} with the use of other variables and Hamiltonian form of equations.  The case of uniformly close orbits of the asteroid and the planet is considered in \cite{lidov-ziglin}. Complete numerical study of bifurcations in the double averaged restricted circular three-body problem is given in  \cite{vashkovyak}. In this problem planar motion is stable with respect to spatial perturbations  {\cite{neish}}. 
 
 Double averaged planar elliptic problem was considered in  \cite{akse} and  \cite{veresh} under assumption  that the distance between  the asteroid and  the star is much larger  than the distance between  the planet and the star.  Complete numerical study of bifurcations in this problem is given in    \cite{vash_planar}. 
 
 These models and results  are of great current interest in  relation to study of motion of exoplanets \cite{shev}.
  
In the restricted elliptic three-body problem there is the family of planar orbits of asteroid, i.e.  orbits  which are in the plane of star-planet system.  For small mass of the planet, majority (in measure sense)  of these orbits  are stable with respect to variations  of initial data in the considered plane. Stability of these orbits  with respect to variations of initial data that put the asteroid out of this plane is an open question.  It can be considered in the framework of the double averaged problem.  In this note we  study  stability of planar orbits which are equilibria of the double averaged  problem. Each such equilibrium is characterised by two parameters: the ratio of the semi-major axes of the asteroid and the planet, and the eccentricity of the orbit of the planet. It is known, that  these orbits are stable  for small enough eccentricities of the planet \cite{neish}. Current interest to this problem is related to the fact, that many exoplanets have large eccentricities and inclinations \cite{skdc,subaru}.

\section{Statement of the problem and the Hamiltonian of the system}

\begin{figure}[htbp] 
\centering
\includegraphics[width=4.0in]{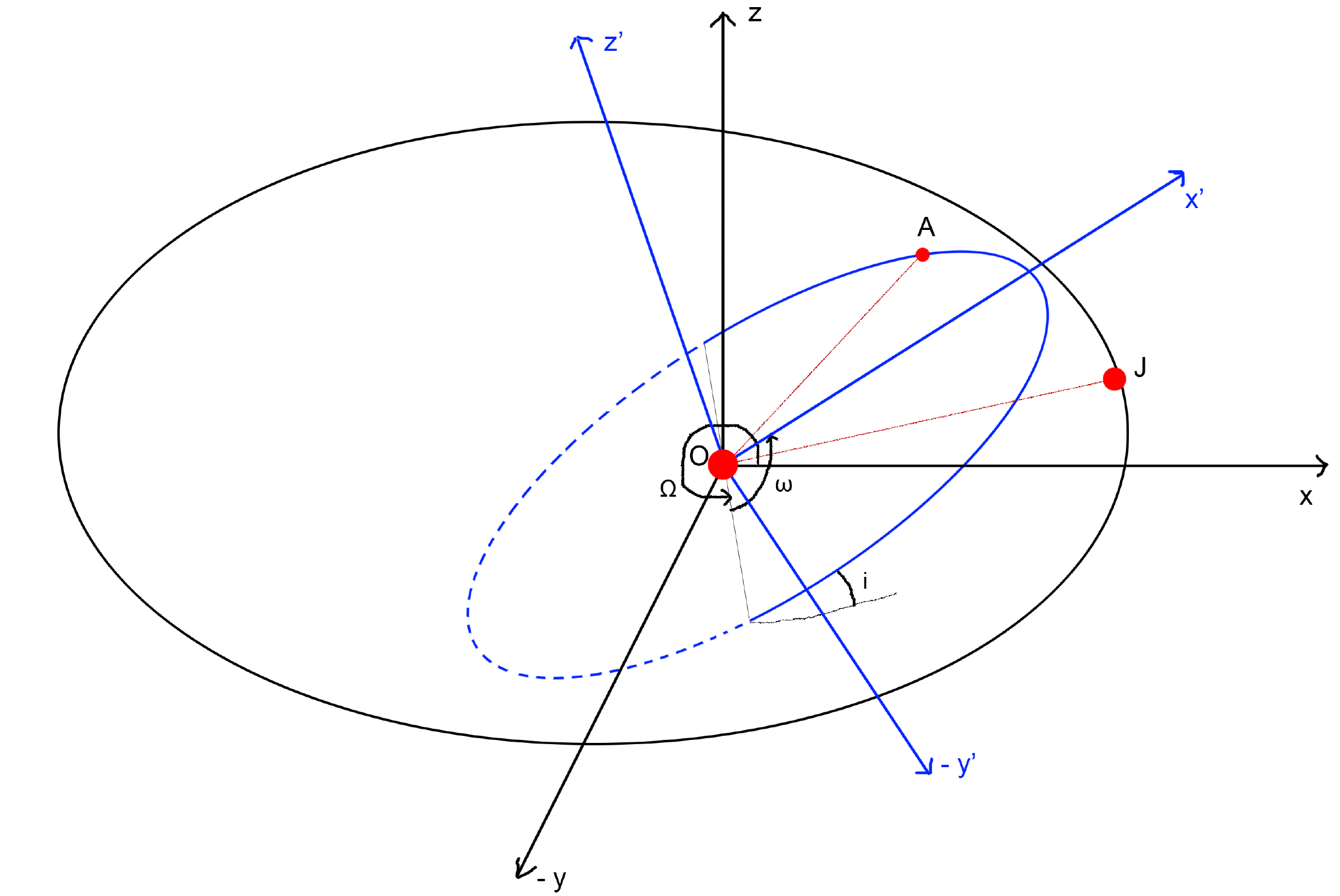} 
\caption{Coordinate frames}
\label{fig:tu77} 
\end{figure} 
Consider spatial restricted elliptic three body problem with a star $S$, planet $J$ and asteroid $A$ \cite{Sze}. Take the origin $O$ of the right Cartesian coordinate system $Oxyz$  at the position of the star  and choose plane of motion of the star  and the planet as the plane $Oxy$ of this system. Let  positive direction of the axis $Ox$ be direction from $O$  towards the periapsis of the orbit of the planet. Let  coordinates of the planet and the asteroid  in this system be   $(x_{J},y_{J},0)$ and $(x,y,z)$, respectively. Denote $a, l, e,\omega, i, \Omega$ the standard osculating elements  of the orbit of the asteroid: the semi-major axis, the mean longitude, the eccentricity, the argument of periapsis, the inclination, the longitude of ascending node.     Introduce a rotating right Cartesian coordinate  frame $Ox'y'z'$ for which the plane   $Ox'y'$ is the osculating  plane of the orbit of the asteroid,  positive direction of the axis $Oz'$ coincides with the direction of the angular momentum of the asteroid about the origin $O$, and positive direction of the axis $Ox'$ is direction from $O$  towards the osculating periapsis of the orbit of the asteroid.   Let $(x',y',0)$ be coordinates of the asteroid in this coordinate frame. Then (see   Fig.\ref{fig:tu77}) 
\begin{equation}
\label{coord_transform}
\left\{
\begin{aligned}
 & x=\left( \cos\Omega
\cos\omega -\cos  i \sin\Omega
  \sin\omega  \right) x'+\left( -\cos\Omega\sin  
\omega  -\cos i \sin\Omega \cos
\omega  \right) y',
\\ 
 & y=\left( \sin\Omega \cos\omega
  +  \cos i\cos\Omega
 \sin\omega  \right) x'+\left( -\sin\Omega\sin
\omega +\cos  i \cos\Omega
  \cos\omega  \right) y' ,
\\ 
 & z=\left(\sin  i
  \sin\omega\right) x' +\left(\sin i \cos\omega\right) y'. 
\end{aligned} 
\right.\
\end{equation}

Take unit of mass such that the sum  of masses of the  star and the planet is $1$. Denote $\mu$   mass of the planet.  The planet moves in a prescribed elliptic orbit:
\begin{eqnarray}
\label{coord_J}
 {{x}_{J}}={{a}_{J}}\left( \cos E_J-{{e}_{J}} \right),\
  {{y}_{J}}={{a}_{J}}\sqrt{1-{{e}_{J}}^{2}}\sin E_J,\quad
  l_J=E_J-e_J\sin E_J. 
\end{eqnarray}
Here $a_J, e_J, E_J, l_J$ are the semi-major axis, the eccentricity,  the eccentric anomaly, and the mean anomaly  of the planet's orbit.  In what follows we put $a_J=1$.

Dynamics of the asteroid can be described using canonical Poincar\'e variables $p_1, p_2, p_3, q_1,q_2,q_3$:  
\begin{equation}
\label{P_elements}
\begin{matrix}
   \begin{aligned}
  & {{p}_{1}}=L ,\\ 
 & {{p}_{2}}=\sqrt{2\left( L-G \right)}\cos\left( g+h \right), \\ 
 & {{p}_{3}}=\sqrt{2\left( G-H \right)}\cos \left( h \right), \\ 
\end{aligned}& 
\begin{aligned}
  & {{q}_{1}}=l+g+h, \\ 
 & {{q}_{2}}=-\sqrt{2\left( L-G \right)}\sin \left( g+h \right), \\ 
 & {{q}_{3}}=-\sqrt{2\left( G-H \right)}\sin \left( h \right), \\ 
\end{aligned}  
\end{matrix}
\end{equation}
where 
$L,G,H, l,g,h$ are canonical Delaynay elements: $L=\sqrt{(1-\mu)a}$, $G=L\sqrt{1-e^{2}}$, $H=G\cos i$,    $l$ is the mean anomaly of asteroid, $g=\omega$, $h=\Omega$
 \cite{Sze}.

The Hamiltonian of the asteroid is  \cite{Sze}
\begin{equation}
\label{F}
\Phi=-{\frac {\left( 1-\mu \right) ^{2}}{2{L}^{2}}}-{\frac {\mu}{\sqrt{(x-x_{J})^{2}+(y-y_{J})^{2}+{z}^{2}}}}-\mu (x\ddot x_J+y\ddot y_J). \end{equation}
Here coordinates $(x,y,z)$ of the asteroid should be expressed via Poincar\'e elements using formulas   (\ref{coord_transform}),   (\ref{P_elements}) and equations of motion of the asteroid in the elliptic orbit:
\begin{eqnarray}
\label{coord_A}
 {{x'}}={{a}}\left( \cos E-{{e}} \right),\
  {{y'}}={{a}}\sqrt{1-{{e}}^{2}}\sin E,\quad
  l=E-e\sin E,
\end{eqnarray}
where $E$ is the eccentric anomaly of the asteroid. Coordinates $x_J, y_J$ of the planet  are prescribed functions of time.  

The double averaged Hamiltonian $\bar \Phi$ is defined as
\begin{equation}
\label{bar_F}
\bar \Phi =\frac{1}{(2\pi)^2}\int_0^{2\pi}\int_0^{2\pi}\Phi\,  {\text d} l{\text d}l_J.
\end{equation}
The double average of the last term in $\Phi$ is 0. Because the double averaged Hamiltonian does not depend on $q_1$, the canonically conjugate variable $p_1=L$ is the first integral of the double averaged system. Thus, the first term in $F$ is constant in  this system. Thus, dynamics of variables $p_1, p_2,  q_1,q_2$ is described by the Hamiltonian system with two degrees of freedom. The Hamiltonian of this system is 
$-\mu \bar V$, where $\mu \bar V$  is double averaged force function  of gravity of the planet:
\begin{equation}
\label{bar_V}
\bar V =\frac{1}{(2\pi)^2}\int_0^{2\pi}\int_0^{2\pi}V{\text d}l {\text d}l_J,\quad V={\frac {1}{\sqrt{(x-x_{J})^{2}+(y-y_{J})^{2}+{z}^{2}}}}.
\end{equation}
The function $V$ depends on two parameters, $a$ and $e_J$ (we take $a_J=1$). 

The double averaged  planar restricted elliptic three-body problem corresponds to the invariant plane $p_3=0, q_3=0$ (i.e. $i=0$) of this problem. Dynamics in this plane is described by the Hamiltonian system with one degree of freedom for the phase variables $p_2, q_2$. Its Hamiltonian depends  on two parameters,   $a$ and $e_J$. Complete numerical study of bifurcations in this problem is given in    \cite{vash_planar}. This system has stable equilibria. In the next Section we discuss stability of these equilibria in the spatial  double averaged   restricted elliptic three-body problem, i.e. stability of these equilibria with respect to spatial perturbations.

\section{Stability of equilibria of double averaged system}
\begin{figure}[htbp] 
\centering
\includegraphics[width=3.0in]{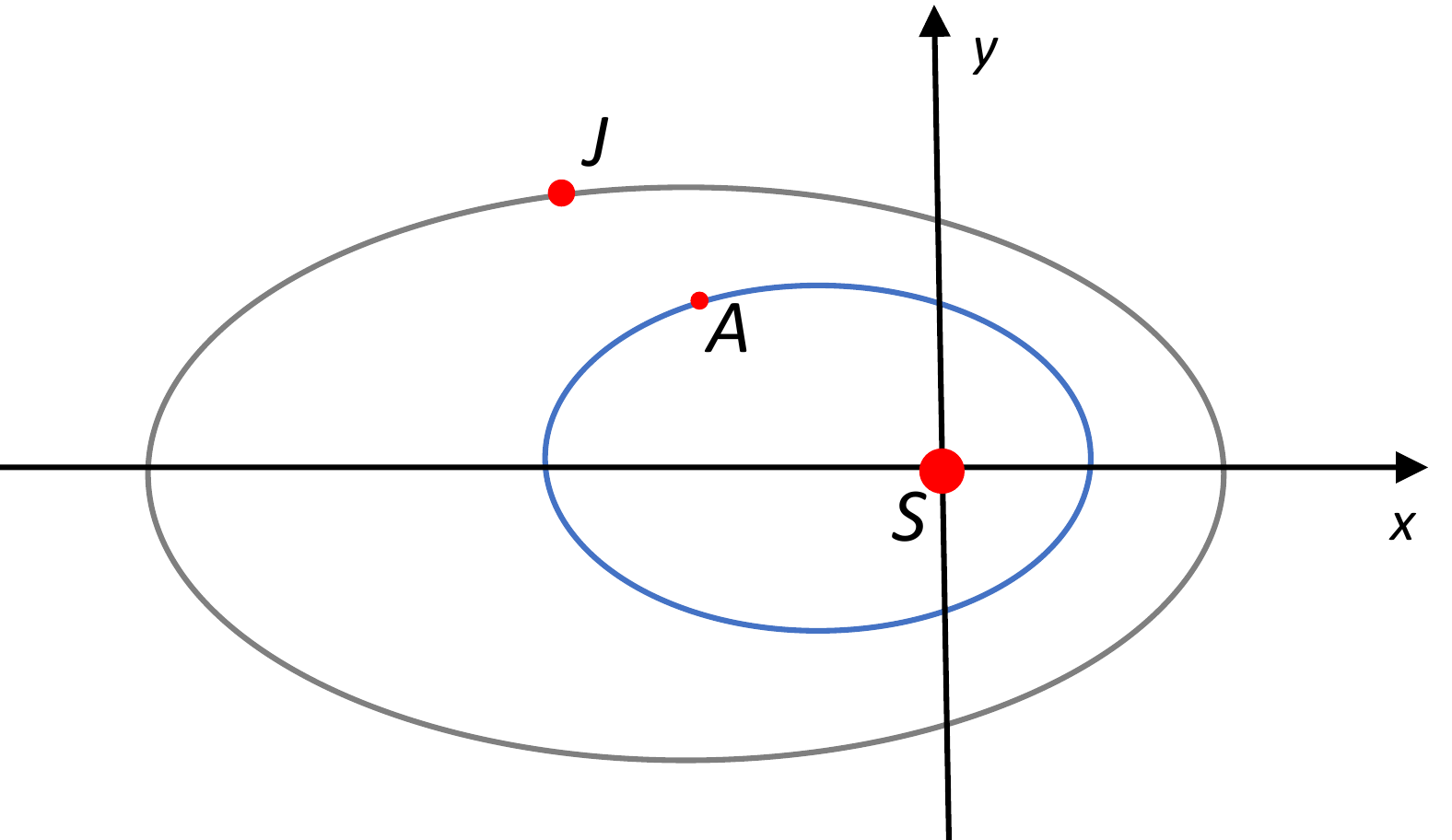} 
\caption{Aligned orbits}
\label{aligned} 
\end{figure}

The double averaged  planar restricted elliptic three-body problem  has stable equilibria for some domains in the plane of parameters $a,e_J$. To study linear stability of these equilibria with respect to spatial perturbations we  consider quadratic  in $p_3,q_3$ part of the function   $\bar V$ at these equilibria. We start with expansion of the function $V$:

\begin{equation*}
V=R+ V_2 +O(p_3^2+ q_3^2), \ V_2= A\, {{p}_{3}}^{2}+2B\, {{p}_{3}}{{q}_{3}}+C\, {{q}_{3}}^{2},
\end{equation*}
where 
\begin{equation}
\begin{aligned}
  & A=-\frac{1}{2}\frac{y{{y}_{J}}}{{{\left( {{\left( x-{{x}_{J}} \right)}^{2}}+{{\left( y-{{y}_{J}} \right)}^{2}} \right)}^{{3}/{2}\;}} G}, \\ 
 & B=-\frac{1}{4}\frac{\left( x{{y}_{J}}+y{{x}_{J}} \right)}{{{\left( {{\left( x-{{x}_{J}} \right)}^{2}}+{{\left( y-{{y}_{J}} \right)}^{2}} \right)}^{{3}/{2}\;}}G}, \\ 
 & C=-\frac{1}{2}\frac{x{{x}_{J}}}{{{\left( {{\left( x-{{x}_{J}} \right)}^{2}}+{{\left( y-{{y}_{J}} \right)}^{2}} \right)}^{{3}/{2}\;}} G}, \\ 
 & R=\frac{1}{\sqrt{{{\left( x-{{x}_{J}} \right)}^{2}}+{{\left( y-{{y}_{J}} \right)}^{2}}}} . 
\end{aligned}
\end{equation}
In calculation of coefficients $A,B,C$ it is taken into account that, as it can be seen  from the phase portraits in \cite{vash_planar},  stable equilibria  of the double averaged planar problem correspond   to $q_{2}=0, p_2>0$, i. e. $g+h=0$, directions from the star to periapses of orbits of the planet and of the asteroid coincide, Fig.\ref{aligned}. These coefficients are calculated at $i=0$.  Thus $x\equiv x', y\equiv y'$ in formulas for these coefficients.

Now we should average $V_2$ over the mean anomaly of the asteroid $l$  and the mean anomaly of the planet $l_J$.

Let us show that the double average of $B$ is 0. Denote $r(l,l_J)=\left(\left( x-{{x}_{J}} \right)^{2}+\left( y-{{y}_{J}} \right)^{2} \right)^{1/2}$. Thus$$
B=-\frac{1}{4G}(B_1+B_2), \ B_1= \frac{ x{{y}_{J}}}{(r(l,l_J))^{3}}, \  B_2= \frac{ y{{x}_{J}}}{(r(l,l_J))^{3}}. 
$$
Symmetry of orbits about $y$-axis implies that   $r(l,l_J)=r(2\pi-l,2\pi-l_J)$, $r(2\pi-l,l_J)=r(l,2\pi-l_J)$.   Then  for the double average of $B_1$ we have
$$
\int_{0}^{2\pi }{}\int_{0}^{2\pi }{B_1}\text{d}l\text{d}{{l}_{J}}=
\int_{0}^{\pi }{}\int_{0}^{\pi } \left(\frac{ x{{y}_{J}}}{(r(l,l_J))^{3}}-\frac{ x{{y}_{J}}}{(r(2\pi-l,2\pi- l_J))^{3}}+\frac{ x{{y}_{J}}}{(r(2\pi-l,l_J))^{3}}-\frac{ x{{y}_{J}}}{(r(l,2\pi-l_J))^{3}}       \right)
\text{d}l\text{d}{{l}_{J}}=0.
$$
Similarly
$$
\int_{0}^{2\pi }{}\int_{0}^{2\pi }{B_2}\text{d}l\text{d}{{l}_{J}}=0.
$$
 Thus the double average of $V_2$ is
\begin{equation}
\bar{V_2}=\bar{A}\,{{p}_{3}}^{2}+\bar{C}\, {{q}_{3}}^{2},
\end{equation}
where 
\begin{equation*}
\begin{aligned}
  & \bar{A}=-\frac{1}{{4{\pi }^{2}G}}
  \int_{0}^{\pi }\int_{0}^{\pi }
  \left( \frac{{{r}_{2}}^{3}-{{r}_{1}}^{3}}{{{r}_{1}}^{3}\, {{r}_{2}}^{3}} \right)y{y}_{J}\text{d}l\text{d}{l}_{J},\\
   & \bar{C}=-\frac{1}{{4{\pi }^{2}G}}
  \int_{0}^{\pi }\int_{0}^{\pi }
  \left( \frac{{{r}_{2}}^{3}+{{r}_{1}}^{3}}{{{r}_{1}}^{3}\, {{r}_{2}}^{3}} \right)x{x}_{J}\text{d}l\text{d}{l}_{J}\\
 \end{aligned}
\end{equation*}
are the average values of $A$ and $C$. 

In above formulas,
\begin{equation*}
\begin{aligned}
  & {{r}_{1}}=\sqrt{{{\left( x-{{x}_{J}} \right)}^{2}}+{{\left( y-{{y}_{J}} \right)}^{2}}}, \\ 
 & {{r}_{2}}=\sqrt{{{\left( x-{{x}_{J}} \right)}^{2}}+{{\left( y+{{y}_{J}} \right)}^{2}}}. \\ 
\end{aligned}
\end{equation*}

For $\bar{A}$, we have 
\begin{equation}
\frac{{{r}_{2}}^{3}-{{r}_{1}}^{3}}{{{r}_{1}}^{3}\, {{r}_{2}}^{3}} >0.
\end{equation}
Thus $\bar{A}$ is negative. 

To study the sign of $\bar C$ we use a numerics. Values $\bar C$ and $\bar R$ (double average of $R$) can be calculated using integration over eccentric anomalies. We have eccentric anomaly $E$ of the asteroid  and eccentric anomaly $E_{J}$ of the planet in our formulas. By Kepler's equation 
\begin{equation}
l=E-e\sin{E}, \ l_J=E_J-e\sin{E_J},
\end{equation}
 for any function $f$ we have
\begin{equation}
\begin{aligned}
   \int_{0}^{2\pi }{}\int_{0}^{2\pi }{f}\text{d}l\text{d}{{l}_{J}}&=\int_{0}^{2\pi }{}\int_{0}^{2\pi }{f}\frac{\text{d}l}{\text{d}E}\text{d}E\frac{\text{d}{{l}_{J}}}{\text{d}{{E}_{J}}}\text{d}{{E}_{J}} \\ 
 & =\int_{0}^{2\pi }{}\int_{0}^{2\pi }{f}\left( 1-e\cos  E  \right)\left( 1-{{e}_{J}}\cos E_J \right)\text{d}E\text{d}{{E}_{J}} . 
\end{aligned}
\end{equation}
Thus we have 
\begin{equation}
\begin{aligned}
 &\bar{C}=-\frac{1}{{4{\pi }^{2}G}}\int_{0}^{\pi }{\int_{0}^{\pi }{\left( \frac{{{r}_{2}}^{3}+{{r}_{1}}^{3}}{{{r}_{1}}^{3}\, {{r}_{2}}^{3}} \right)x{{x}_{J}} \left( 1-e\cos  E  \right)\left( 1-{{e}_{J}}\cos E_J \right) }}\text{d}E\text{d}{{E}_{J}} ,\\ 
 & \bar{R}=\frac{1}{{2{\pi }^{2}}}\int_{0}^{\pi }{\int_{0}^{\pi }{\left( \frac{{{r}_{1}}+{{r}_{2}}}{{{r}_{1}} {{r}_{2}}} \right)}{{ \left( 1-e\cos  E  \right)\left( 1-{{e}_{J}}\cos E_J \right) }}}\text{d}E\text{d}{{E}_{J}}. 
\end{aligned}
\end{equation}
Note that $\mu \bar{R}$ is the double averaged force function of gravity of the planet in the planar elliptic three-body problem. Equilibria of the planar problem are points in the plane $p_2, q_2$ such that $q_2=0$ and ${\partial \bar{R}}/{\partial p_2}=0$. This is equivalent to
 ${\partial \bar{R}}/{\partial e}=0$ at $g=0$.  We  calculated $\bar{C}$ numerically  when $a_{J}=1$,  $a$, $e_{J}$  are on some grid, and $e$ is found from the equation ${\partial \bar{R}}/{\partial e}=0$. Numerically, we found that $\bar{C}$ is always negative. This was checked analytically in the limit cases of a small $a$ and of orbits close to collision orbits.

Thus, $\bar A<0$,  $\bar C <0$, and  $\bar V_2$ is  a negative definite quadratic form.  Hence, stable equilibria of the double averaged  planar restricted elliptic three-body problem are stable in the linear approximation  as equilibria of the  double averaged spatial restricted elliptic three-body problem for all values of parameters. It follows from results of \cite{akse}, \cite{vash_planar} that the quadratic form of expansion of $\bar R$ near  stable equilibria of  the double averaged  planar restricted elliptic three-body problem  is positive definite for all values of parameters. Therefore, quadratic terms of our expansion do not provide a  Lyapunov function for study stability. However, Arnold-Moser theorem in Kolmogorov-Arnold-Moser (KAM)  theory guarantee Lyapunov stability of equilibria of systems with two degrees of freedom if a) there are no resonances between frequencies up to 4th order, and b) some non-degeneracy property for  4th order terms in the normal form near the equilibrium is satisfied (\cite{akn}, Sec. 8.3.3).  Using expansion for small $a$ one can check that for small $a$ conditions a) and b) are violated on  some curves in parameters plane only.  Because our system is analytic, this implies that  conditions a) and b) can be violated on some curves in the parameters plane only. Thus, the considered equilibria are  Lyapunov stable for all values of parameters $a,e_J$ except, possibly,  for parameters belonging to some finite  set of analytic curves.

\section{Conclusion}

We have shown that stable equilibria of the double averaged  planar restricted elliptic three-body problem are linearly stable as equilibria of the double averaged  spatial  restricted elliptic three-body problem.   KAM theory implies that these equilibria are Lyapunov stable for all parameters but, possibly, parameters from some finite set of analytic curves.  These exceptional values  of parameters correspond to a finite set of resonances and to a degeneration.   In a separate note we will show that indeed there is an instability for a resonance 2:1 between frequencies of oscillations in the plane of star-planet system and across this plane.


\end{document}